\documentclass{article}
\usepackage{amsmath,amssymb,amsthm,dsfont,hyperref,tikz}
\begin{document}
	\title{Lie Conformal Algebra and Dual Pair Type Realizations of Some Moonshine Type VOAs, and Calculations of the Correlation Functions}
	\author{Hongbo Zhao}
	\date{}
	\maketitle
	\allowdisplaybreaks
	\newtheorem{theorem}{Theorem}
	\newtheorem{definition}{Definition}
	\newtheorem{example}{Example}
	\newtheorem{lemma}{Lemma}
	\newtheorem{proposition}{Proposition}
	\numberwithin{equation}{section}
\section{Introduction}

\par{
	A vertex operator algebra (VOA) $V$ is called a CFT type VOA if $V$ is $\mathbb{Z}_{\geq 0}$-graded with $V_0=\mathbb{C}1$. A moonshine type VOA is a CFT type VOA such that $V_1=\{0\}$. It is well known that for a moonshine type VOA $V$, the space $V_2$ is has a commutative but not necessarily associative algebraic structure with the product given by
	$$
		a\circ b:=a(1)b,
	$$
	and we call $V_2$ the Greiss algebra of $V$.
}
\par{
	Finite dimensional Jordan algebras are examples of algebras which are commutative but in general not assocaitive, so we are interested in moonshine type VOAs whose Greiss algebras are finite dimensional simple Jordan algebras. In particular, we focus on Hermitian type Jordan algebras.
}
\par{
	The examples of moonshine type VOAs whose Greiss algebras are Hermitian type Jordan algebras were first given by Lam in \cite{Lam1} and \cite{Lam2}. In \cite{AM}, Ashihara and Miyamoto constructed VOAs $V_{\mathcal{J},r}$ parametrized by $r\in\mathbb{C}$ whose Greiss algebras are type $B$ Jordan algebras $\mathcal{J}$ and properties of $V_{\mathcal{J},r}$ were studied by Niibori and Sagaki in \cite{NS}.
	In \cite{Z2}, we constructed the simple VOAs $\bar{V}_{\mathcal{J},r}, r\in\mathbb{Z}_{\neq 0}$ for type $B$ Jordan algebras, and we constructed $V_{\mathcal{J},r}$ for the remaining cases of Hermitian Jordan algebras in \cite{Z3}.
}
\par{
	This paper is a continuation of the papers \cite{Z1}, \cite{Z2}, \cite{Z3}. The main tools we use here are Lie conformal algebras (LCA) and certain free fields. We give LCA realizations of the VOAs $V_{\mathcal{J},r}$, and free field realizations of the simple VOAs $\bar{V}_{\mathcal{J},r}$ using dual pair type constructions. Moreover, we give a differenet approach to calculating correlation functions of the VOAs $V_{\mathcal{J},r}$, using our realization of the simple VOAs $\bar{V}_{\mathcal{J},r}$.
}
\par{
	Recall that a Lie conformal algebra (LCA) is a $\mathbb{C}[\partial]$-module $V$ equipped with a $\mathbb{C}$-linear map called the $\lambda$-bracket $[\cdot_{\lambda}\cdot]:V\otimes V\rightarrow \mathbb{C}[\lambda]\otimes V$ such that for any $a,b,c\in V$ the following identities hold:
	\begin{align*}
	&[(\partial a)_{\lambda}b]=-\lambda[a_{\lambda}b],\,[a_{\lambda}b]=-[b_{-\lambda-\partial}a],\\
	&[a_{\lambda}[b_{\mu}c]]=[b_{\mu}[a_{\lambda}c]]+[[a_{\lambda}b]_{\lambda+\mu}c].
	\end{align*}
	It is noted that we can get an LCA from a vertex algebra by forgetting the negative integer products and only consider the non-negative integer products. More precisely, let $V$ be a vertex algebra, then $V$ is also an LCA with $\partial=T$ and the $\lambda$-bracket
	\begin{align}
	b_{\lambda}a=\sum_{k\geq 0}\frac{\lambda^k}{k!}b(k)a\label{LCA}.
	\end{align}
}
\par{
	For our purpose we need the concept of vertex operator super algebras (VOSAs) which generalizes VOAs slightly, and fields we need are also understood in the super sense. Following \cite{Kac}, given a space $V$, a collection of fields $S(z)=\{a(z)|\,a(z)\in End(V)[[z,z^{-1}]]\}$ over $V$ is called a free field, if for any three fields $a(z),b(z),c(z)\in S(z)$ we have
	$$
		[[a(z),b(w)],c(x)]=0.
	$$
	Free fields generate many important VO(S)As. The most fundamental examples are $bc$-system VOSAs and $\beta\gamma$-system VOAs. For our purpose, the free fields we need are the generating fields of Heisenberg VOAs and Symplectic Fermions VOSAs. They are analogues of the free fields generating $bc$-system and $\beta\gamma$-system VOSAs.
}
\par{
	The main content of this paper is organized as follows. In Section 2, we introduce some LCAs together with their central extensions, and in Section 3 we construct VOAs $V_{\mathcal{J},r}$ using the LCAs introduced in Section 2, and the first main result is Theorem 1. In all these constructions we use free fields, and the way we construct $V_{\mathcal{J},r}$ is different from the ones given in \cite{AM} and \cite{Z2}. In Section 4 we constrcut the simple VOAs $\bar{V}_{\mathcal{J},r}$ where $r\in \mathbb{Z}_{\neq 0}$, and they are all dual pair type constructions, which are summarized in Theorem 2 and Theorem 3. As an application, we calculate the correlation functions of generating fields of the VOAs $V_{\mathcal{J},r}$, giving Theorem 4, where free fields and the dual pair constructions play key roles.
}
\section{Heisenberg and Symplectic Fermion VO(S)As, and Some Lie Conformal Algebras}
\par{
	In this section we briefly review the construction of certain free fields which in particular generate Heisenberg and Symplectic Fermion VO(S)As. To be convenient we use the language of superspaces. Throughout this paper, the ground field is $\mathbb{C}$ and all concepts involving vector spaces are understood in the super sense.
}
\par{
	Let $\mathfrak{h}=\mathfrak{h}_{\bar{0}}\oplus \mathfrak{h}_{\bar{1}}$ be a superspace with
	$$
		sdim(\mathfrak{h})=(p|q),
	$$
	such that $\mathfrak{h}$ is equipped with a non-degenerate supersymmetric bilienar form $(\cdot,\cdot)$. We define the Lie superalgebra
	$$
		\hat{\mathfrak{h}}:=\mathfrak{h}\otimes \mathbb{C}[t,t^{-1}]\oplus\mathbb{C}K
	$$
	with the Lie superbracket
	$$
		[a(m),b(n)]=m\delta_{m+n,0}(a,b)K,\,[a(m),K]=0,\,a(m):=a\otimes t^m.
	$$
	The even part of this Lie superalgebra is
	$$
		\mathfrak{h}_{\bar{0}}\otimes \mathbb{C}[t,t^{-1}]\oplus\mathbb{C}K
	$$
	and the odd part is
	$$
		\mathfrak{h}_{\bar{1}}\otimes \mathbb{C}[t,t^{-1}].
	$$
	The sub Lie superalgebra
	$$
		\hat{\mathfrak{h}}_-:=\mathfrak{h}\otimes t^{-1}\mathbb{C}[t^{-1}]\oplus\mathbb{C}K
	$$
	is supercommutative, and we check that the supersymmetric tensor space
	$$
		S(\hat{\mathfrak{h}}_-)
	$$
	is a VOSA. The generaing fields are
	$$
		Y(a,z)=\sum_{k\in\mathbb{Z}}a(k)z^{-k-1},
	$$
	and the whole VOSA structure is obtained from the normal ordered products. It is clear that 
	$$
		\{a(z)|\,a\in\mathfrak{h}\}
	$$
	forms a free field, and this will be used in Section 5.
}
\par{
	When $\mathfrak{h}$ is purely even the corresponding VOA is exactly the Heisenberg VOA, while $\mathfrak{h}$ is purely odd the corresponding VOSA is the symplectic Fermion VOSA. As we mentioned in the introduction, there are another two important examples called $bc$-systems and $\beta\gamma$ systems, but we do not need these two examples in this paper and their relations will be discussed elsewhere.
}
\par{
	Recall that for any vertex algebra $V$, it is automatically an LCA by (\ref{LCA}). We note that for the VOSAs $S(\hat{\mathfrak{h}}_-)$, the subspace
	$$
		\widehat{S_2(\mathfrak{h})}:=span\{a(-m)b(-n)\cdot 1,1|\,a,b\in V\}
	$$
	is even which is also closed under the $k$-th products for all $k\geq 0$, therefore $\widehat{S_2(\mathfrak{h})}$ is an LCA. The one dimensional subspace $\mathbb{C}1$ is central, and the quotient space
	$$
		S_2(\mathfrak{h})=\widehat{S_2(\mathfrak{h})}/\mathbb{C}1\simeq span\{a(-m)b(-n)\cdot 1|\,a,b\in V\}
	$$
	is also an LCA. Moreover, $\widehat{S_2(\mathfrak{h})}$ is a nontrivial central extension of $S_2(\mathfrak{h})$, and for $x,y\in S_2(\mathfrak{h})$ the corresponding non-trivial $2$-cocycle $c(x,y)$ is given by
	\begin{align}
		c(x,y):=\frac{1}{(|x|+|y|-1)!}\hat{x}(|x|+|y|-1)\hat{y}\label{cocycle}
	\end{align}
	where $\hat{x}$ and $\hat{y}$ denote arbitrary choices of preimages of $x$ and $y$ in $\widehat{S_2(\mathfrak{h})}$, and for $x=a(-m)b(-n)\cdot 1$
	$$
		|x|:=m+n
	$$
	denotes the degree of $x$. In other words, let 
	$$
		\pi :\widehat{S_2(\mathfrak{h})}\rightarrow \mathbb{C}1
	$$
	denote the projection onto the space $\mathbb{C}1$, then
	$$
		c(x,y)\lambda^{|x|+|y|-1}=\pi([\hat{x}_{\lambda}\hat{y}]).
	$$
}
\par{
	For our purpose we also need to introduce another LCA describing the type $A$ case. Let $\mathfrak{h}$ be an even vector space and $\mathfrak{h}^*$ be the corresponding dual space, then there is a unique way to assign a supersymmetric bilinear form $(\cdot,\cdot)$ on $\mathfrak{h}\oplus \mathfrak{h}^*$ such that $\mathfrak{h},\mathfrak{h}^*$ are isotropic subspaces, $\mathfrak{h}\oplus \mathfrak{h}^*$ is even, and
	$$
		(a,b^*)=\langle b^*,a\rangle.
	$$
	We note that there is a $\mathbb{C}^{\times}$-action on $\mathfrak{h}\oplus \mathfrak{h}^*$:
	\begin{align}
		c\cdot(a,b^*)=(ca,c^{-1}b^*)\label{caction},
	\end{align}
	and we have two fixed point LCAs $S_2(\mathfrak{h}\oplus \mathfrak{h}^*)^{\mathbb{C}^{\times}}$ and $\widehat{S_2(\mathfrak{h}\oplus \mathfrak{h}^*)}^{\mathbb{C}^{\times}}$.
	It follows that 
	\begin{align*}
	&S_2(\mathfrak{h}\oplus \mathfrak{h}^*)^{\mathbb{C}^{\times}}=span\{a(-m)b^*(-n)\cdot 1|\, a\in \mathfrak{h},\,b^*\in \mathfrak{h}^* \},\\
	&\widehat{S_2(\mathfrak{h}\oplus \mathfrak{h}^*)^{\mathbb{C}^{\times}}}:=span\{a(-m)b^*(-n)\cdot 1,1|\, a\in \mathfrak{h},\,b^*\in \mathfrak{h}^*\}.
	\end{align*}
	We remark that we can also impose $\mathfrak{h}$ and $\mathfrak{h}\oplus \mathfrak{h}^*$ to be odd spaces, but there is only a minus sign difference on the level and we omit this equivalent choice.
}
\par{
	We also remark that the LCAs $S_2(\mathfrak{h})$ and $S_2(\mathfrak{h}\oplus \mathfrak{h}^*)^{\mathbb{C}^{\times}}$  introduced here are not new. It can be shown that if $sdim(\mathfrak{h})=(N|0)$ or $(0|N)$, then
	\begin{align}
		S_2(\mathfrak{h}\oplus \mathfrak{h}^*)^{\mathbb{C}^{\times}}\simeq gc_{N,Ix},\label{gcn}
	\end{align}
	and
	\begin{align}
		S_2(\mathfrak{h})\simeq \begin{cases}
			oc_{N,Jx},\,\mathfrak{h}=(0|2N),\\
			spc_{N,Ix},\,\mathfrak{h}=(N|0),
		\end{cases}\label{on}
	\end{align}
	where $gc_{N,Ix}$, $oc_{N,Jx}$ and $spc_{N,Ix}$ are certain sub LCAs of the general linear conformal algebra $gc_N$ introduced in \cite{BKL03}. Moreover, all these LCAs are simple, and the $2$-cocycle $c(x,y)$ is the same as the $2$-cocyle introduced in \cite{Kac} up to a constant scalar. 
}
\par{
	For $oc_{N,jx}$ and $spc_{N,Ix}$ we set
	$$
		L_{a,b}(m,n)1：=\frac{1}{2}a(-m)b(-n)\cdot 1,
	$$
	and in particular
	$$
		L_{a,b}:=L_{a,b}(-1,-1)1.
	$$
	Similarly for $gc_{N,Ix}$, 	
	$$
		L_{a,b^*}(m,n)1：=\frac{1}{2}a(-m)b^*(-n)\cdot 1,\,
		L_{a,b^*}:=L_{a,b^*}(-1,-1)1.
	$$
	The following lemma will be used later
	\begin{lemma}
	For $N\geq 2$, the LCAs
	$gc_{N,Ix}$, $oc_{N,Jx}$ and $spc_{N,Ix}$ are generated by their degree two subspaces
	\begin{align*}
		&S_2(\mathfrak{h})_2:=span\{L_{a,b}|a,b\in \mathfrak{h}\}\,\\&S_2(\mathfrak{h}\oplus \mathfrak{h}^*)^{\mathbb{C}^{\times}}_2:=span\{L_{a,b^*}|a\in \mathfrak{h},\,b^*\in \mathfrak{h}^*\},
	\end{align*}
	through the identifications (\ref{gcn}), (\ref{on}), as LCAs.
	\end{lemma}
}
\par{
	\textbf{Proof.} We only prove the case for $spc_{N,Ix}$ with $N=2$, $\mathfrak{h}$ purely even, and general cases are proved in a similar way. The computation is similar to the proof of Theorem 1.1, (3), of \cite{Z1}. We assume that $\mathfrak{h}$ is spanned by the basis $a,b$ with
	$$
		(a,a)=(b,b)=1,\,(a,b)=0.
	$$
	We note that $L_{a,a},L_{b,b}$ are mutually orthogonal Virasoro elements, and
	$$
		L_{a,b}(-m,-n)1=\frac{L_{a,a}(0)^{m-1}L_{b,b}(0)^{n-1}}{(m-1)!(n-1)!}L_{a,b},	
	$$
	therefore $L_{a,b}(-m,-n)1$ can be generated. We also note that
        \begin{align*}
	&L^r_{b,a}(1)L^r_{a,b}(-m,-n)1
	=mL^r_{b,b}(-m,-n) 1+nL^r_{a,a}(-m,-n)1,\\
	&(L^r_{b,a}(-2,-1)1)(2)L^r_{a,b}(-m,-n)1
	=m(m-1)L^r_{b,b}(-m,-n)\c1-n(n+1)L^r_{a,a}(-m,-n)1
\end{align*}
hold by a direct computation, therefore $L_{a,a}(-m,-n)1$ and $L_{b,b}(-m,-n)1$ can also be generated by solving the equation, and we conclude the proof.
}
\par{
	We have mentioned in \cite{Z2} and $\cite{Z3}$ that the condition $N\geq 2$ is necessary. 
}
\section{Constructing $V_{\mathcal{J},r}$ Using Lie Conformal Algebras}
\par{
	In this section we use Lie conformal algebras $gc_{N,Ix}$, $oc_{N,Jx}$,  $spc_{N,Ix}$, and the 2-cocycle (\ref{cocycle}) to construct the VOAs $V_{\mathcal{J},r}$.
}
\par{
	We first briefly review some facts about Hermitian Jordan algebras which will be used later. Let $\mathcal{J}_X$ be a Hermitian Jordan algebra where $X$ denotes the type, $X=A,B,C$. Let $d$ denote the rank of $\mathcal{J}_X$, then $\mathcal{J}_X$ can be described using tensors. Let $\mathfrak{h}$ be a superspace with a non-degenerate supersymmetric bilinear form $(\cdot,\cdot)$ according to the type $X$:
	\begin{align}
	sdim(\mathfrak{h})=\begin{cases}
	(d|0),\,X=A,\\
	(0|2d),\,X=C,\\
	(d|0),\,X=B.
	\end{cases}\label{type}
	\end{align}
	Then we identify the dual space $\mathfrak{h}^*$ with $\mathfrak{h}$, and as vector spaces the Hermitian Jordan algebras are
	$$
	\mathcal{J}_X=\begin{cases}
	\mathfrak{h}\otimes \mathfrak{h}^*\simeq S^2(\mathfrak{h}\oplus\mathfrak{h}^*)^{\mathbb{C}^{\times}},\,X=A,\\
	S^2(\mathfrak{h}),\,X=B,C.
	\end{cases}
	$$
	We introduce elements
	\begin{align}
		L_{a,b^*}:=a\otimes b^*\in \mathcal{J}_A \label{typea}
	\end{align}
	in type $A$ case, and 
	\begin{align}
		L_{a,b}:=\begin{cases}
		a\otimes b+b\otimes a,\,X=B,\\
		a\otimes b-b\otimes a,\,X=C
		\end{cases}\label{typebc}
	\end{align}
	in $B,C$ cases. Then
	$$
		\mathcal{J}_X=\begin{cases}
		 span\{L_{a,b^*}|\,a\in\mathfrak{h},b^*\in\mathfrak{h}^*\},\,X=A,\\
		span\{L_{a,b}|\,a,b\in \mathfrak{h}\},\,X=B,C.
	\end{cases}
	$$
	and in any cases the Jordan product is given by
	$$
	A\circ B:=\frac{1}{2}(AB+BA)
	$$
	where $A,B$ are viewed as elements in $End(\mathfrak{h})\simeq \mathfrak{h}\otimes \mathfrak{h}$ with the associative product
	$$
	(a\otimes b)(u\otimes v)=(b,u)(a\otimes v).
	$$
}
\par{
	We briefly recall the construction of $V_{\mathcal{J}_X,r}$ , where the case for $X=B$ was given in \cite{AM} and the cases for $X=B,C$ were given in \cite{Z3}. The details can be found in \cite{Z3}, Section 3 and Section 4. Let $\mathcal{L}_X$ be the infinite dimensional Lie algebra which is chosen according to the type $X$:
	$$
	\mathcal{L}_X=\begin{cases}
	\mathfrak{gl}_{\infty},\,X=A,\\\mathfrak{sp}_{\infty},\,X=B,\\
	\mathfrak{so}_{\infty},\,X=C.\\
	\end{cases}	
	$$
	We write $\mathcal{L}$ for short to agree with the notation in \cite{Z3} if there is no ambiguity. By making a certain Lie algebra decomposition
	\begin{align*}
	&\mathcal{L}=\mathcal{L}_{-}\bigoplus\mathcal{L}_{+}
	\end{align*}
	and defining a one dimensional $\mathcal{L}_+$-module spanned by $v_r$, we have an induced $\mathcal{L}$-module $M_r$:
	\begin{align*}
	M_r:=&U(\mathcal{L})\otimes_{U(\mathcal{L}_{+})}\mathbb{C}v_r
	\cong U(\mathcal{L}_{-})v_r.
	\end{align*}
	It was shown in \cite{AM} and \cite{Z3} that $M_r$ is a vertex operator algebra, denoted by $V_{\mathcal{J}_X,r}$.
}
\par{
	Now we give the LCA realization of $V_{\mathcal{J},r}$, which is the main result of this section. We first review a general construction of vertex algebras starting from certain LCAs with a non-trivial 2-cocyle. Recall that we can construct `formal distribution Lie algebras', and moreover, vertex algebras from LCAs (\cite{Kac}, Section 2.7). Let $C$ e an LCA. We consider the `affinization'
	$$
	\hat{C}:=C\otimes \mathbb{C}[t,t^{-1}],
	$$
	and define the quotient space
	$$
	Lie(C):=\hat{C}/span\{(\partial a)t^n+nat^{n-1}|\,a\in C\}.
	$$
	Then $Lie(C)$ is a Lie algebra with the following Lie bracket which can be verified from the axioms of LCAs:
	$$
	[at^m,bt^n]:=\sum_{k\geq 0}{m\choose k}(a(k)b)t^{m+n-k}.
	$$
	Moreover, if $C$ is a free $\mathbb{C}[\partial]$-module generated by $A$,
	$$
	C=\mathbb{C}[\partial]\otimes A,
	$$
	then
	$$
	Lie(C)\simeq A\otimes \mathbb{C}[t,t^{-1}]
	$$
	with a derivation $\partial=-\partial_t$.
	For each $a\in A$ we can associate a $Lie(C)$-valued formal power series
	$$
	a(z):=\sum_{k}(at^k)z^{-k-1},
	$$
	and it is checked that any two power series $a(z)$ and $b(w)$ are mutually local.
}
\par{
	Let $c(\cdot,\cdot)$ be a $2$-cocyle of $C$ and $K$ be the central element. In our cases, it is always possitble to find a non-trivial one with $\partial K=0$. Then 
	$$
	Kt^n=0\text{ for all }n\neq -1,
	$$
	and we have a centrally extended Lie algebra 
	$$
	Lie(C\bigoplus\mathbb{C}K)=Lie(C)\oplus\mathbb{C}Kt^{-1}
	$$
	with the new Lie bracket $[\cdot,\cdot]'$:
	$$
		[x,y]'=[x,y]+c(x,y)Kt^{-1}.
	$$ 
	Take a split
	$$
	Lie(C\bigoplus\mathbb{C}K):=Lie(C\bigoplus\mathbb{C}K)_+\oplus Lie(C\bigoplus\mathbb{C}K)_-
	$$
	where
	\begin{align*}
	&Lie(C\bigoplus\mathbb{C}K)_+:=span\{at^n|\,a\in A,\,n\geq 0\}\oplus\mathbb{C}Kt^{-1},\\
	&Lie(C\bigoplus\mathbb{C}K)_-:=span\{at^n|\,a\in A,\,n< 0\}.
	\end{align*}
	Define a one dimensional $Lie(C\bigoplus\mathbb{C}K)_+$-module spanned by $1$:
	$$
	at^n\cdot 1=0,\,\text{ for all }a\in AK,\,n\geq 0,\,Kt^{-1}\cdot 1=r\cdot 1.
	$$
	Then we construct a (level $r$) vertex algebra as an induced $Lie(C\bigoplus\mathbb{C}K)$-module
	$$
	Vert_r(C):=U(Lie(C\bigoplus\mathbb{C}K))\otimes_{U(Lie(C\bigoplus\mathbb{C}K)_+)}\mathbb{C}1\simeq U(Lie(C\bigoplus\mathbb{C}K)_-)\cdot 1
	$$
	as we note that the $Lie(C)$-valued power series $a(z)$, $a\in C$ are generating mutually local fields over $Vert_r(C)$. The whole vertex algebra structure on $Vert_r(C)$ is determined by normal ordered products according the the reconstruction theorem \cite{Kac}. It follows that $C$ is a sub space of $Vert_r(C)$ through
	$$
		C\rightarrow Ct^{-1}\subseteq Vert_r(C),
	$$
	which will be used in Section 5.
}
\par{
	For a CFT type VOA $V$, we can get an LCA $V/\mathbb{C}1$ and a corresponding 2-cocyle $c(\cdot,\cdot)$ in the same way as (\ref{cocycle}). We note that in general, the vertex algebra
	$$
		Vert_r(V/\mathbb{C}1)
	$$
	is very `large' and not very interesting. A modification is that we take a reasonable sub LCA of $V$ instead, and for particular examples of Heisenberg and symplectic Fermion VO(S)As, we choose sub LCAs $gc_{N,Ix}$, $oc_{N,Jx}$ and $spc_{N,Ix}$.
}

\par{
	We choose the LCA according to the type $X$. Let $C_X$ denote the following LCAs according to the type $X$:
	$$
		C_X=\begin{cases}
		S_2(\mathfrak{h}),\,X=B,C,\\
		S_2(\mathfrak{h}\oplus\mathfrak{h}^*)^{\mathbb{C}^{\times}},\,X=A.
		\end{cases}
	$$
	in which the type of $\mathfrak{h}$ is determined by (\ref{type}). Let $c(\cdot,\cdot)$ be the 2-cocycle (\ref{cocycle}), then we have the vertex algebra
	$$
		Vert_r(C_X).
	$$
}
\par{
	We have the following lemma
	\begin{lemma}
		The Lie algebra $Lie(C_X\oplus \mathbb{C}K)$ acts on $M_r$ with $K=r$
	\end{lemma}	
	This is verified by a direct comparison and calculation.
}
\par{
	Using Lemma 2 the following theorem is a direct consequence of our construction.
	\begin{theorem}
		The $Lie(C_X\oplus \mathbb{C}K)$-action on $M_r$ extends to a VOA isomorphism
		$$
			Vert_r(C_X)\rightarrow V_{\mathcal{J}_X,r}.
		$$
		which is induced by
		$$
			\begin{cases}
				L_{a,b}t^{-1}\cdot 1\mapsto L_{a,b},\,X=B,C,\\
				L_{a,b^*}t^{-1}\cdot 1\mapsto L_{a,b^*},\,X=A.
			\end{cases}
		$$
	\end{theorem}
}
\par{
	The following proposition follows easily from Lemma 1 and Theorem 1, and the case for $X=B$ was proved in \cite{NS} as Proposition 3.1.
	\begin{proposition}
		The VOAs $V=V_{\mathcal{J}_X,r}$ are generated by their degree two subspaces $V_2$.
	\end{proposition}
}
\section{Dual Pair Construction of the Simple VOAs $\bar{V}_{\mathcal{J},r}$}
\par{
	In \cite{NS} the authors proved that the VOA $V_{\mathcal{J}_B,r}$ is simple if and only if $r\notin \mathbb{Z}$, and in \cite{Z2}, we give the explicit construction of the simple quotients $\bar{V}_{\mathcal{J}_B,r}$ for all $r\in\mathbb{Z}_{\neq 0}$ using dual pair constructions. In \cite{Z3} we show that $V_{\mathcal{J}_A,r}$ and $V_{\mathcal{J}_C,r}$ are also simple if and only if $r\notin \mathbb{Z}$. The main result of this section is to give a uniform description of the simple quotients $\bar{V}_{\mathcal{J}_X,r}$, $X=A,B,C$.
}
\par{
	The main idea is the following: Let $\mathfrak{h}$ and $\mathfrak{h}'$ be two superspaces of types $(p|q)$ and $(p'|q')$ , with non-degenerate supersymmetric bilinear forms $(\cdot,\cdot)$ and $(\cdot,\cdot)'$ respectively. Then the tensor product space
	$$
		\mathfrak{h}\otimes\mathfrak{h}'
	$$
	is a supersymmetric space of type $(pp'+qq'|pq'+p'q)$ with the non-degenerate supersymmetric bilinear form $(\cdot,\cdot)\otimes(\cdot,\cdot)'$. On the level of VOSAs, we use the VOSA associated with $\mathfrak{h}\otimes \mathfrak{h}'$ described in Section 2, and we consider the action of the `Lie supergroup' $Osp(\mathfrak{h}')$ on $\mathfrak{h}'$, where the Lie supergroup $Osp(\mathfrak{h}')$ is understood as a Harish-Chandra pair \cite{Z2}. To construct the simple quotients of the VOA $V_{\mathcal{J}_X,r}$ when $r\in \mathbb{Z}_{\neq 0}$, the type of $\mathfrak{h}$ is related to the type $X$ of the Jordan algebra, and the type of $\mathfrak{h}'$ is related to the level $r$.
}
\par{
	We first describe the cases $X=B,C$. Let $d$ denote the rank of the Jordan algebra $\mathcal{J}_X$, and $r$ be the level, then we choose the type of $\mathfrak{h}$ and $\mathfrak{h}'$ as follows
	\begin{align*}
		sdim(\mathfrak{h})=
		\begin{cases}
			(d|0),\,X=B,\\
			(0|2d),\,X=C,
		\end{cases}
	\end{align*}
	and
	\begin{align*}
		sdim(\mathfrak{h}')=
		\begin{cases}
			(r|0),\,r\in\mathbb{Z}_{\geq 1},\\
			(0|-r),\,r<0,\,r\in 2\mathbb{Z},\\
			(1|-r+1),\,r<0,\,r\in 2\mathbb{Z}+1.
		\end{cases}
	\end{align*}
	It is clear that for our choices of $\mathfrak{h}'$, the following relation always holds
	$$
		r=p'-q'.
	$$
}
\begin{theorem}
	For $X=B,C$ and $r\in \mathbb{Z}_{\neq 0}$,
	$$	
		\bar{V}_{\mathcal{J}_X,r}\simeq S(\widehat{(\mathfrak{h}\otimes \mathfrak{h}')}_-)^{Osp(\mathfrak{h}')}.
	$$
\end{theorem}
\par{
	We describe $\bar{V}_{\mathcal{J}_X,r}$ more explicitly using invariant theory for Lie supergroups $Osp(\mathfrak{h}')$ (See for example, \cite{LZ1},\cite{Serg}). For any $m,n\in\mathbb{Z}$ we introduce elements
	$$
		L^r_{a,b}(m,n):=\sum_i(a\otimes e_i)(m)(b\otimes e_i^*)(n)
	$$
	which are elements in the enveloping algebra $U((\widehat{\mathfrak{h}\otimes \mathfrak{h}'})_-)$. It follows that 
	\begin{align*}
		 &S(\widehat{(\mathfrak{h}\otimes \mathfrak{h}')}_-)^{Osp(\mathfrak{h}')}\\
		=&{\rm span}\{L^{r}_{a_1,b_1}(-m_1,-n_1)\cdots L^{r}_{a_k,b_k}(-m_k,-n_k)\cdot 1|\,a_i,b_i\in\mathfrak{h},m_i,n_i\geq 1\},
	\end{align*}
	and we also use $L^r_{a,b}$ to denote the element
	$$
	L^r_{a,b}(-1,-1)\cdot1\in
	S(\widehat{(\mathfrak{h}\otimes \mathfrak{h}')}_-)^{Osp(\mathfrak{h}')}
	$$
	if there is no confusion.
}
\par{
	\begin{proposition}
		For $X=B,C$, there is a VOA homomorphism
		$$
			V_{\mathcal{J}_X,r}\rightarrow S(\widehat{(\mathfrak{h}\otimes \mathfrak{h}')}_-)^{Osp(\mathfrak{h}')}.
		$$
		which is induced by
		$$
			\phi:L_{a,b}\mapsto L^r_{a,b}.
		$$
	\end{proposition}
}
\par{
	\textbf{Proof.} It is easy to verify that for all $k\geq 0$ we have
	$$
	\phi(L_{a,b}(k)L_{u,v})=\phi(L_{a,b})(k)\phi(L_{u,v})
	$$
	by setting $K=r$, therefore it extends to an LCA homomorphism
	$$
		\phi:L_{a,b}(-m,-n)1\mapsto L^r_{a,b}(-m,-n)1,\,1\mapsto 1
	$$
	and moreover a VOA homomorphism 
	$$
	V_{\mathcal{J}_X,r}\rightarrow S(\widehat{(\mathfrak{h}\otimes \mathfrak{h}')}_-)^{Osp(\mathfrak{h}')}
	$$
	by Theorem 1 and Proposition 1.
}
\par{
	It is clear that $V_{\mathcal{J},r}$ has a unique simple quotient, and to conclude the proof of Theorem 2 we need to prove the following.
	\begin{proposition}
		The fixed point VOA $$	
			 S(\widehat{(\mathfrak{h}\otimes \mathfrak{h}')}_-)^{Osp(\mathfrak{h}')}.
		$$
		is simple.
	\end{proposition}
}
\par{
	The proof for type $B$ case has already been given in \cite{Z2}, and the proof of the general case is similar. First, we note that the VO(S)A
	$$
		S(\widehat{(\mathfrak{h}\otimes \mathfrak{h}')}_-)
	$$
	is simple, and there is an invariant bilienar form over them \cite{KR87}. The main point is that in any cases, $Osp(\mathfrak{h}')$ acts semisimplely on $S(\widehat{(\mathfrak{h}\otimes \mathfrak{h}')}_-)$ which also preserves the invariant bilinear form. For simplicity we set
	$$
		M:=S(\widehat{(\mathfrak{h}\otimes \mathfrak{h}')}_-),
	$$
	and there is a decomposition
	$$
		M=M^{Osp(\mathfrak{h}')}\oplus M'
	$$
	where
	$$
		M'=\oplus_{\lambda}M^{\lambda}
	$$
	is the direct sum of non-trivial irreducible $Osp(\mathfrak{h}')$-modules. Therefore, $M^{Osp(\mathfrak{h}')}$ is orthogonal to $M'$. The non-degeneracy of the invariant bilinear form on $M$ implies that it is also non-degenerate on $M^{Osp(\mathfrak{h}')}$, and therefore the fixed point VOA is simple by \cite{Li}.
}
\par{
	For type $A$ case the result is similar with a slight difference.
\begin{theorem}
	For $r\in \mathbb{Z}_{\neq 0}$,
	$$	
	\bar{V}_{\mathcal{J}_A,r}\simeq S(\widehat{(\mathfrak{h}\oplus \mathfrak{h}^*)\otimes \mathfrak{h}')}_-)^{Osp(\mathfrak{h}')\times \mathbb{C}^{\times}}.
	$$
	Here the $\mathbb{C}^{\times}$-action is the one induced by (\ref{caction}).
\end{theorem}
This theorem is proved by the corresponding analogues of Proposition 2 and Proposition 3. The proof is similar by considering the $\mathbb{C}^{\times}$-action on $\mathfrak{h}\oplus \mathfrak{h}^*$, and we omit the details.
}
\section{Calculation of the Correlation Functions Using $\bar{V}_{\mathcal{J},r}$}
\par{
	In this section we apply the simple VOAs $\bar{V}_{\mathcal{J}_X,r}$, $r\in\mathbb{Z}_{\neq 0}$ to compute the correlation functions of the generating fields in the VOAs $V_{\mathcal{J}_X,r}$.
}
\par{
	We introduce some notations and conventions. For a VOA $V$, $v\in V$, we write
	$$
		v(z):=Y(a,z)
	$$
	for short, and suppose the mode expasion is 
	$$
		v(z)=\sum_k v(k)z^{-k-1},
	$$
	we set
	$$
		v(z)_+:=\sum_{k\geq 0} v(k)z^{-k-1},\,v(z)_-:=\sum_{k< 0} v(k)z^{-k-1}.
	$$
	In this section, all VO(S)As $V$ are assumed to be of CFT type, and we use $1'$ to denote the unique element in the restricted dual $V^*$ such that
	$$
		\langle 1',1\rangle=1
	$$
	and
	$$
		\langle 1',v\rangle =0
	$$
	for all $v\in V_{i}$, $i>0$. 
}
\par{
	For $v_1,\cdots,v_n\in V$, we call
	\begin{align}
	\langle 1',v_1(z_1)\cdots v_{n}(z_{n})1\rangle\label{corr}
	\end{align}
	the correlation function (or genus zero $n$-point function) of the fields $v_1(z_1),\cdots,v_n(z_n)$. It can be viewed either as an element in $\mathbb{C}((z_1))\cdots((z_n))$, or as a rational function of $(z_i-z_j),i\neq j$. As a formal power series can be obtained by the corresponding multivariable function on the domain
	$$
		|z_1|>\cdots|z_n|.
	$$
	It is also well known that for any permutation $\tau\in S_n$ (\ref{corr}) equals
	$$
		\langle 1',v_{\tau(1)}(z_{\tau(1)})\cdots v_{\tau(n)}(z_{\tau(n)})1\rangle
	$$
	as multivariable complex functions after doing analytic continuations.
}
\par{
	The main result of this section is to compute (\ref{corr}) for
	$$
		v_i=\begin{cases}
			v_i=L_{a_i,b_i},X=B,C,\\
			v_i=L_{a_i,b^*_i},X=A.
		\end{cases}
	$$
	The type $B$ case has been computed in \cite{Z1}, but the method we use here is different.
}
\par{
	We first recall the following recursion formula for CFT type VOAs, which is a variant of Lemma 2.2.1 in \cite{Z}
	\begin{lemma}
		Suppose $v_i\in V_k$, $k\geq 1$, and $n\geq 2$, then
	\begin{align}
		&\langle 1',v_1(z_1)\cdots v_{n}(z_{n})1\rangle\notag\\
		=&\sum_{j=2,\cdots, n,k\geq 0}\iota_{z_1,z_j}(z_1-z_j)^{-k-1}\langle 1',v_2(z_2)\cdots (v_1(k)v_j)(z_j)\cdots v_{n}(z_{n})1\rangle.\label{recur}
	\end{align}
	\end{lemma}
	\textbf{Proof.} We note that because $V$ is of CFT type,
	$$
	\langle 1',v_1(z_1)_-\cdots v_{n}(z_{n})1\rangle=0
	$$
	by counting the degree, hence
	$$
	\langle 1',v_1(z_1)\cdots v_{n}(z_{n})1\rangle=\langle 1',v_1(z_1)_+\cdots v_{n}(z_{n})1\rangle.
	$$
	We also recall that
	$$
	[v_s(z)_+,v_t(w)]=\sum_{j\geq 0}(v_s(j)v_t)(w)\iota_{z,w}(z-w)^{-j-1},
	$$
	and
	$$
	v(z)_+1=0.
	$$
    Therefore the recursion formula by obtained by moving $v_1(z_1)_+$ to the rightmost position.
}
\par{
	\begin{lemma}
		As a complex function the correlation function
		$$
			\langle 1', L_{a_1,b_1}(z_1)\cdots L_{a_n,b_n}(z_n)1\rangle
		$$
		has the form
		$$
			\mathbb{C}[r][(z_i-z_j)^{-1}|\,i<j]
		$$
		and the constant coefficients only depend on $a_i,b_i$. Moreover as a polynomial in $r$
		the degree of the correlation function is less than $n$.
	\end{lemma}
	}
\par{
	To prove the lemma we need to prove a slightly stronger version which allows that $v_i$ are elements in $C_X\subseteq Vert_r(C_X)\simeq V_{\mathcal{J}_X,r}$, as $C_X$ is closed under the $\lambda$-brackets: Suppose $v_i=L_{a_i,b_i}(-m_i,-n_i)1$, then the conclusion of Lemma 4 still holds. We prove this by doing induction on $n$. When $n=0,1,2$ the conclusion holds by a direct computation. Suppose the conclusion holds for $n=k$. then we note that for all $l\geq 0$ and $i,j=1,\cdots,n$
	$$
		v_i(l)v_j=v+rv'
	$$
	for some $v,v'\in C_X$ which can be determined, therefore we conclude the proof by the induction hypothesis and formula (\ref{recur}).
}
\par{
	As another corollary of the recursion formula, we have
	\begin{lemma}
		Suppose There is a VOA homomorphism $V\rightarrow V'$, $v_i\in V$, and $v'_i$ are the homomorphic images of $v_i$ in $V$, then the correlation function
		$$
		\langle 1',v_1(z_1)\cdots v_n(z_n)1\rangle
		$$
		equals the correlation function
		$$
		\langle 1',v'_1(z_1)\cdots v'_n(z_n)1\rangle.
		$$
	\end{lemma}
}
	\par{
		The recursion formula in the proof of Lemma 4 cannot be used directly for the computation and we solve this in another way. Recall that by Proposition 2 and Proposition 3 there are homomorphisms
		$$
			V_{\mathcal{J}_X,r}\rightarrow \bar{V}_{\mathcal{J}_X,r}
		$$
		and when $r\in\mathbb{Z}_{\neq 0}$ the simple VOAs $\bar{V}_{\mathcal{J}_X,r}$ are realized
		using free fields and dual pair type constructions, therefore by Lemma 6,
		$$
			\langle 1', L_{a_1,b_1}(z_1)\cdots L_{a_n,b_n}(z_n)1\rangle
			=\langle 1', L^r_{a_1,b_1}(z_1)\cdots L^r_{a_n,b_n}(z_n)1\rangle
		$$
		for $r\in\mathbb{Z}_{\geq 0}$. Because a complex coefficient polynomial of degree less than $n$ is uniquely determined by its values
		at $n$ different points and we can view the correlation functions as complex functions, therefore by Lemma 5 we only need to calculate the correlation functions $\bar{V}_{\mathcal{J}_X,r}$,
		for $n$ different non-zero integer value of $r$. In particular it is sufficient to compute
		$$
			\langle 1', L^r_{a_1,b_1}(z_1)\cdots L^r_{a_n,b_n}(z_n)1\rangle
		$$
		for all positive integers $r\in\mathbb{Z}_{>0}$.
	}
	\par{
		We need the following Theorem for free fields.
		\begin{proposition}[Wick's Theorem, \cite{Kac}, Theorem 3.3]
			Let $a_1(z),\cdots,a_m(z)$ and $b_1(w),\cdots,b_n(w)$ be a collections of
			free fields. Then
			\begin{align*}
				&:a_1(z)\cdots a_m(z)::b_1(w)\cdots b_n(w):\\
				=&\sum_{s=0}^{min(m,n)}\sum_{(i_1,j_1),\cdots,(i_s,j_s)}(-1)^{\epsilon(i_1,j_1;\cdots;i_s,j_s)}[a_{i_1}(z)_+,b_{j_1}(w)_-]\cdots [a_{i_s}(z)_+,b_{j_s}(w)_-]\\
				&:a_1(z)\cdots a_m(z)b_1(w)\cdots b_n(w):_{(i_1,j_1;\cdots;i_s,j_s)}.
			\end{align*}
			Here the subscript $(i_1,j_1;\cdots;i_s,j_s)$ in the last line means that the fields $a_{i_1}(z),b_{j_1}(w),\\\cdots a_{i_s}(z),b_{j_s}(w)$ are removed, and $\epsilon(\cdot)=\pm 1$ is the sign obtained by the super rule: each permutation of the adjacent odd fields changes the sign.
		\end{proposition}
	}
	\par{
		For two free field $a(z)$ and $b(w)$ we call
		$$
			[a(z)_+,b(w)_-]
		$$
		the (Wick's) contraction of $a(z)$ and $b(w)$. In particular, for free fields $a(z)$, $b(w)$ in the Heisenberg and symplectic Fermion VOSAs, it is clear that
		$$
			[a(z)_+,b(w)]=[a(z)_+,b(w)_-]=\iota_{z,w}\frac{(a,b)}{(z-w)^2}.
		$$
	}
	\par{
		The following lemma is obtained by counting the degree of a CFT type VOA.
		\begin{lemma}
			$$
			\langle 1',:u_1(z_1)\cdots u_m(z_m):1\rangle=0.
			$$
		\end{lemma}
	}
	\par{
		Now we can use Wick's theorem to compute the correlation function 
		$$
		\langle 1', L^r_{a_1,b_1}(z_1)\cdots L^r_{a_n,b_n}(z_n)1\rangle.
		$$
		Because 
		\begin{align*}
			L^r_{a,b}(z)=&\frac{1}{2}\sum_i ((a\otimes e_i)(-1)(b\otimes e_i))(z)\\
			=&\frac{1}{2}\sum_i :(a\otimes e_i)(z)(b\otimes e_i)(z):
		\end{align*}
		by definition, we can use Proposition 4 repeatedly. 
	}
	\par{
		Before doing the computation we briefly recall some notations in \cite{Z1} and we generalize them slightly. For a sequence of fields $L_{a_1,b_1}(z_1),\cdots,L_{a_n,b_n}(z_1)$ in the cases $X=B,C$, we have a corresponding sequence, denoted by
		$$
			T=(a_1,b_1)\cdots(a_n,b_n).
		$$
		A $BC$-type diagram over the sequence $T=(a_1,b_1)\cdots(a_n,b_n)$ is a graph, with the vertex set $V=\{a_1,b_1,\cdots,a_n,b_n\}$, and edge set $E$ consisting of unordered pairs $\{u,v\},\;u,v\in V$ satisfying:
        \begin{description}
            \item[(1).] $\{a_i,b_i\},\{a_i,a_i\},\{b_i,b_i\}\notin E$ for all $i=1,\cdots,n$.
            \item[(2).] Any two edges have no common point.
            \item[(3).] $|E|=n$.
        \end{description}
		Denote the set of all $BC$-type diagrams over $T$ by $D(T)$. Similarly in the $A$-type case, a sequence of fields $L_{a_1,b^*_1}(z_1),\cdots,L_{a_n,b^*_n}(z_1)$ gives rise to a sequence
		$$
			T=(a_1,b^*_1)\cdots(a_n,b^*_n)
		$$
		and a diagram over $T$
		is a $BC$-type diagram over $T$ satisfying one more condition:
		\begin{description}
            \item[(1)*.] $\{a_i,a_j\},\{b^*_i,b^*_j\}\notin E$ for all $1\leq i<j\leq n$.
		\end{description}
		We also denote the set of $A$-type diagram over $T$ by $D(T)$ if there is no ambiguity. 
}
\par{
	As an example, for $T=(a_1,b_1)(a_2,b_2)(a_3,b_3)(a_4,b_4)$, the following
	$$
	
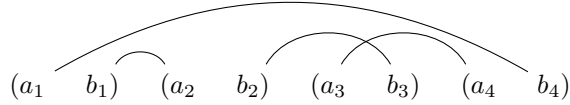
\begin{tikzpicture}
	\node (p1) at (-3,0) {$(a_1$};
	\node (p2) at (-2,0) {$b_1)$};
	\node (p3) at (-1,0) {$(a_2$};
	\node (p4) at (0,0) {$b_2)$};
	\node (p5) at (1,0) {$(a_3$};
	\node (p6) at (2,0) {$b_3)$};
	\node (p7) at (3,0) {$(a_4$};
	\node (p8) at (4,0) {$b_4)$};
	\draw (p2)to [out=60,in=120](p3);
	\draw (p1)to [out=30,in=150](p8);
	\draw (p4)to [out=60,in=120](p6);
	\draw (p5)to [out=60,in=120](p7);
	\end{tikzpicture}
	$$
	is a $BC$-type diagram (over $T$), and it corresponds to the following contraction which gives the term
	\begin{align}
		&\frac{1}{2^4}\sum_{i_1,i_2,i_3,i_4}[(a_1\otimes e_{i_1})(z_1)_+,(b_4\otimes e_{i_4})(z_4)_-][(b_1\otimes e_{i_1})(z_1)_+,(a_2\otimes e_{i_2})(z_2)_-]\notag\\&[(b_2\otimes e_{i_2})(z_2)_+,(b_3\otimes e_{i_3})(z_3)_-][(a_3\otimes e_{i_3})(z_3)_+,(a_4\otimes e_{i_4})(z_4)_-]\langle 1',1\rangle \notag\\
		=&\frac{r(a_1,b_4)(a_2,b_1)(b_2,b_3)(a_3,a_4)}{16(z_1-z_4)^2(z_1-z_2)^2(z_2-z_3)^2(z_3-z_4)^2}.\label{comp}
	\end{align}
	On the other hand, the following diagram
	$$
	
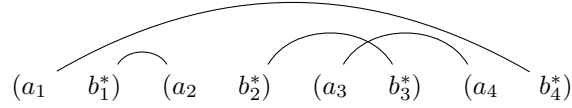
\begin{tikzpicture}
	\node (p1) at (-3,0) {$(a_1$};
	\node (p2) at (-2,0) {$b^*_1)$};
	\node (p3) at (-1,0) {$(a_2$};
	\node (p4) at (0,0) {$b^*_2)$};
	\node (p5) at (1,0) {$(a_3$};
	\node (p6) at (2,0) {$b^*_3)$};
	\node (p7) at (3,0) {$(a_4$};
	\node (p8) at (4,0) {$b^*_4)$};
	\draw (p2)to [out=60,in=120](p3);
	\draw (p1)to [out=30,in=150](p8);
	\draw (p4)to [out=60,in=120](p6);
	\draw (p5)to [out=60,in=120](p7);
	\end{tikzpicture}
	$$
	is not a $A$-type diagram over $T=(a_1,b^*_1)(a_2,b^*_2)(a_3,b^*_3)(a_4,b^*_4)$, but
	$$
	
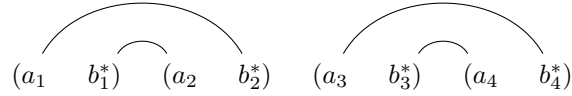
\begin{tikzpicture}
	\node (p1) at (-3,0) {$(a_1$};
	\node (p2) at (-2,0) {$b_1^*)$};
	\node (p3) at (-1,0) {$(a_2$};
	\node (p4) at (0,0) {$b_2^*)$};
	\node (p5) at (1,0) {$(a_3$};
	\node (p6) at (2,0) {$b_3^*)$};
	\node (p7) at (3,0) {$(a_4$};
	\node (p8) at (4,0) {$b_4^*)$};
	\draw (p1)to [out=60,in=120](p4);
	\draw (p2)to [out=60,in=120](p3);
	\draw (p5)to [out=60,in=120](p8);
	\draw (p6)to [out=60,in=120](p7);
	\end{tikzpicture}
	$$
	is $A$-type diagram over $T=(a_1,b^*_1)(a_2,b^*_2)(a_3,b^*_3)(a_4,b^*_4)$.
}
\par{
	We note that for each diagram $D$,
		We have a corresponding $n$-vertex graph $\sigma_D$ by collapsing $a_i,b_i$
		to a single vertex labelled by $i$. More precisely, $\sigma_D$ has vertex set $\{1,\cdots,n\}$ and if $(a_i,a_j)$, $(a_i,b_j)$,
		or $(b_i,b_j)$ is an edge of $D$, then $(i,j)$ is an edge of $\sigma_D$. It is clear that
		 $\sigma_D$ is a disjoint union of cycles.
		}
	\par{
		By Lemma 7, it is sufficient to consider the terms in the summation such that 
		all free fields are contracted, and each of these terms corresponds exactly to a diagram
		$D\in D(T)$. More precisely, $a(z_i)$ is contracted with $b(z_j)$ if and only if $e:=(a,b)$ is an
		edge of $D$, and in this case we define
		$$
			K(e,Z):=\frac{1}{(z_i-z_j)^2}.
		$$
		We also define
		$$
			\Gamma(D):=\sum_{\{a,b\}\in E}(a,b),
		$$
		then by a direct computation similar to (\ref{comp}), a single contraction term corresponding to a diagram $D$ equals
		$$
			\Gamma(D)r^{c(\sigma_D)}\prod_{e\in E}K(e,Z),
		$$
		where $c(\cdot)$ denotes the number of cycles. 
		}
	\par{
		Therefore by taking the summation over $D$ we have
		\begin{align*}
			&		\langle 1', L^r_{a_1,b_1}(z_1)\cdots L^r_{a_n,b_n}(z_n)1\rangle\\
			=&\sum_{D\in D(T)}\Gamma(D)r^{c(\sigma_D)}\prod_{e\in E}K(e,Z)\\
		&\sum_{\sigma\in DR(T)}\Gamma(\sigma,T)\Gamma(\sigma;Z)r^{c(\sigma)},\label{form1}
		\end{align*}
		Here $DR(T)$ denote the set of permutations in $n$ elements such that each block has size $\geq 2$, and $\Gamma(\sigma,T)$, $\Gamma(\sigma,Z)$ means the following: Assume that $\sigma=(C_1)\cdots(C_s)=(k_{11}\cdots k_{1 t_1})\cdots (k_{s1}\cdots k_{s t_s})$, then
\begin{align*}
	\Gamma_{BC}(\sigma,T)\stackrel{def}{=}&2^{-s-n}\prod_{i=1}^{s}Tr(L_{a_{k_{i1}},b_{k_{i1}}}\cdots L_{a_{k_{it_i}},b_{k_{it_i}}}),\\
	\Gamma(\sigma;Z)\stackrel{def}{=}&\prod_{i=1}^n\frac{1}{(z_i-z_{\sigma(i)})^2},
\end{align*}
where $L_{a,b}$ are understood as elements in $End(\mathfrak{h})$ in the sense of (\ref{typebc}). This result coincides with \cite{Z1}, Theorem 1.
We note that in the type $C$ case, although we need to take into account of the sign,
the formula is the same.
}
\par{
	The formula for the type $A$ case is slightly different. For $\sigma=(C_1)\cdots(C_s)=(k_{11}\cdots k_{1 t_1})\cdots (k_{s1}\cdots k_{s t_s})$ we define
	\begin{align*}
	\Gamma_A(\sigma,T)\stackrel{def}{=}2^{-n}\prod_{i=1}^{s}Tr(L_{a_{k_{i1}},b^*_{k_{i1}}}\cdots L_{a_{k_{it_i}},b^*_{k_{it_i}}}).
	\end{align*}
	where $L_{a,b^*}$ are understood as elements in $End(\mathfrak{h})$ in the sense of (\ref{typea}), and we just need to note that a term corresponds to non-type A diagram $D$ equals zero. Then we have
	\begin{align*}
	&\langle 1', L_{a_1,b^*_1}(z_1)\cdots L_{a_n,b^*_n}(z_n)1\rangle\\
	&\sum_{\sigma\in DR(T)}\Gamma_A(\sigma,T)\Gamma(\sigma;Z)r^{c(\sigma)}.
	\end{align*}
}
\par{
	In Summary we have
	\begin{theorem}
		The correlation function
		$
		\langle 1', L_{a_1,b_1}(z_1)\cdots L_{a_n,b_n}(z_n)1\rangle
		$
		in $BC$ type cases or
		$
		\langle 1', L_{a_1,b^*_1}(z_1)\cdots L_{a_n,b^*_n}(z_n)1\rangle
		$
		in A type case equals
		\begin{align*}
		\sum_{\sigma\in DR(T)}\Gamma_X(\sigma,T)\Gamma(\sigma;Z)r^{c(\sigma)}.
		\end{align*}
	\end{theorem}
}

\end{document}